\documentclass[12pt,leqno]{amsart}

\usepackage{amsmath, amsthm, amscd, amsfonts, amssymb, graphicx}
\usepackage[bookmarksnumbered, colorlinks, plainpages]{hyperref}

\numberwithin{equation}{section}

\newtheorem{thm}{Theorem}[section]

\newtheorem{prop}[thm]{Proposition}
\newtheorem{cor}[thm]{Corollary}
\theoremstyle{definition}

\newtheorem{prob}[thm]{Problem}

\newtheorem{rem}[thm]{Remark}
\newtheorem*{leA}{Lemma A}

\def\d{\delta}

\def\<{\langle}
\def\>{\rangle}

\newcommand{\Norm}[1]{\left\|{#1}\right\|}

\def\1{\mathop{\mathrm{(i)}}\nolimits}
\def\2{\mathop{\mathrm{(ii)}}\nolimits}
\def\3{\mathop{\mathrm{(iii)}}\nolimits}

\allowdisplaybreaks

\begin{document}

%%%%%%%%%%%%%%%%%%%%%%%%%%%%%%%%%%%%%%%%%%%%%%%%%%%%%%%%%           M A K E T I T L E                             %%%%%%%%%%%%%%%%%%%%%%%%%%%%%%%%%%%%%%%%%%%%%%%%%%%%%%%%
\title[Operator parallelogram law]{Schatten $p$-norm inequalities related to an extended operator parallelogram law}

\author[M.S. Moslehian, M. Tominaga, K.-S. Saito]{Mohammad Sal Moslehian$^1$, Masaru Tominaga$^2$, Kichi-Suke Saito$^3$}

\address{$^1$ Department of Mathematics, Center of Excellence in Analysis on Algebraic Structures (CEAAS), Ferdowsi University of Mashhad, P.O. Box 1159, Mashhad 91775, Iran.}
\email{moslehian@ferdowsi.um.ac.ir and moslehian@ams.org}
\address{$^2$ Department of Information Systems and Management, Faculty of Applied Information Science, Hiroshima Institute of Technology, 2-1-1 Miyake, Saeki-ku, Hiroshima 731-5193, Japan.}
\email{m.tominaga.3n@it-hiroshima.ac.jp}

\address{$^3$ Department of Mathematics, Faculty of Science, Niigata University, Niigata 950-2181, Japan.}
\email{saito@math.sc.niigata-u.ac.jp}



\maketitle

%%%%%%%%%%%%%%%%%%%%%%%%%%%%%%%%%%%%%%%%%%%%%%%%%%%%%%%%%           A B S T R A C T                               %%%%%%%%%%%%%%%%%%%%%%%%%%%%%%%%%%%%%%%%%%%%%%%%%%%%%%%%
\begin{abstract}
Let $\mathcal{C}_p$ be the Schatten $p$-class for $p>0$.
Generalizations of the parallelogram law for the Schatten $2$-norms have been given in the following form:
If $\mathbf{A}=\{A_1,A_2,\ldots,A_n\}$ and $\mathbf{B}=\{B_1,B_2,\ldots,B_n\}$ are two sets of operators in $\mathcal{C}_2$, then
$$\sum_{i,j=1}^n\|A_i-A_j\|_2^2
+ \sum_{i,j=1}^n\|B_i-B_j\|_2^2
= 2\sum_{i,j=1}^n\|A_i-B_j\|_2^2
- 2\Norm{\sum_{i=1}^n(A_i-B_i)}_2^2.$$
In this paper, we give generalizations of this as pairs of inequalities for Schatten $p$-norms,
which hold for certain values of $p$ and reduce to the equality
above for $p=2$. Moreover, we present some related inequalities for three sets of operators.
\end{abstract}
%%%%%%%%%%%%%%%%%%%%%%%%%%%%%%%%%%%%%%%%%%%%%%%%%%%%%%%%

%%%%%%%%%%%%%%%%%%%%%%%%%%%%%%%%%%%%%%%%%%%%%%%%%%%%%%%%%          1.  I N T R O D U C T I O N                    %%%%%%%%%%%%%%%%%%%%%%%%%%%%%%%%%%%%%%%%%%%%%%%%%%%%%%%%
\section{Introduction}

Suppose that ${\mathbb B}({\mathcal H})$ denotes the algebra of all bounded linear operators on a separable complex Hilbert space ${\mathcal H}$ endowed with an inner product $\langle \cdot , \cdot \rangle$. Let $A\in \mathbb{B}(\mathcal{H})$ be a compact operator and let $\{s_j(A)\}$ denote the sequence of decreasingly ordered singular values of $A$, i.e. the eigenvalues of $|A|=(A^*A)^{1/2}$. The Schatten $p$-norm ($p$-quasi-norm, resp.) for $1 \le p < \infty$ ($0<p<1$, resp.) is defined by $$\|A\|_p=\left(\sum_{j=1}^\infty s_j^p(A)\right)^{1/p}.$$
For $p>0$, the Schatten $p$-class, denoted by $\mathcal{C}_p$, is defined to be the two-sided ideal in $\mathbb{B}(\mathcal{H})$ of those compact operators $A$ for which $\|A\|_p$ is finite.
Clearly
\begin{eqnarray}\label{mos}
\left\|\,|A|^2\,\right\|_{p/2}=\|A\|_p^2
\end{eqnarray}
 for $p>0$.
In particular, $\mathcal{C}_1$ and $\mathcal{C}_2$ are the trace class and the Hilbert-Schmidt class, respectively.
For $1 \le p < \infty$, $\mathcal{C}_p$ is a Banach space; in particular the triangle inequality holds.
For $0<p<1$, the quasi-norm $\|\cdot\|_p$ does not satisfy the triangle inequality, but however satisfies the inequality $\|A+B\|_p^p \leq \|A\|_p^p+\|B\|_p^p$.
For more information on the theory of Schatten $p$-norms the reader is referred to \cite[Chapter 2]{79.BS}.

It follows from \cite[Corollary 2.7]{10.MM} that for $A_1, \ldots, A_n, B_1, \ldots, B_n \in {\mathbb B}({\mathcal H})$
\begin{eqnarray} \label{ineq:norm}
\quad
\sum_{i,j=1}^n\|A_i-A_j\|^2
+ \sum_{i,j=1}^n\|B_i-B_j\|^2
= 2\sum_{i,j=1}^n\|A_i-B_j\|^2
- 2\Norm{\sum_{i=1}^n(A_i-B_i)}^2,
\end{eqnarray}
which is indeed a generalization of the classical \emph{parallelogram law}:
$$|z+w|^2+|z-w|^2=2|z|^2+2|w|^2 \qquad (z,w \in \mathbb{C}).$$
There are several extensions of parallelogram law among them we could refer the interested reader to \cite{03.TC.VC.MT.JP.79, 95.OE.1, 86.MK.37, 82.WZ.23, HZ.MF}.
Generalizations of the parallelogram law for the Schatten $p$-norms have been given in the form of the celebrated Clarkson inequalities (see \cite{08.OH.FK.369} and references therein).
Since $\mathcal{C}_2$ is a Hilbert space under the inner product $\langle A, B\rangle=\textrm{tr} (B^*A)$, it follows from an equality similar to \eqref{ineq:norm} stated for vectors of a Hilbert space (see \cite[Corollary 2.7]{10.MM}) that if $A_1, \ldots, A_n, B_1, \ldots, B_n \in \mathcal{C}_2$, then
\begin{equation} \label{ineq:2-norm}
\sum_{i,j=1}^n\|A_i-A_j\|_2^2
+ \sum_{i,j=1}^n\|B_i-B_j\|_2^2
= 2\sum_{i,j=1}^n\|A_i-B_j\|_2^2
- 2\Norm{\sum_{i=1}^n(A_i-B_i)}_2^2.
\end{equation}
In \cite{10.MM} a joint operator extension of the Bohr and parallelogram inequalities is presented. In particular, it follows from  \cite[Corollary 2.3]{10.MM} that if $A_1, \ldots, A_n, B_1, \ldots, B_n \in \mathbb{B}(\mathcal{H})$, then
\begin{equation} \label{oit}
\begin{split}
\sum_{1 \le i<j \le n} \left| A_i-A_j \right|^2
+
\sum_{1 \le i<j \le n} \left| B_i-B_j \right|^2 =
\sum_{i,j=1}^n \left| A_i - B_j\right|^2
- \left| \sum_{i=1}^n(A_i-B_i) \right|^2.
\end{split}
\end{equation}
In this paper, we give a generalization of the equality \eqref{ineq:2-norm} for the Schatten $p$-norms ($p>0$). First we present similar consideration for three sets of operators. In addition, we provide pairs of complementary inequalities that reduce to \eqref{ineq:2-norm} for the certain value $p=2$.

%%%%%%%%%%%%%%%%%%%%%%%%%%%%%%%%%%%%%%%%%%%%%%%%%%%%%%%%%          2.  SCHATTEN $p$-NORM INEQUALITIES             %%%%%%%%%%%%%%%%%%%%%%%%%%%%%%%%%%%%%%%%%%%%%%%%%%%%%%%%
\section{Schatten $p$-norm inequalities}

To accomplish our results, we need the following lemma which can be deduced from \cite[Lemma 4]{90.RB.FK.719} and \cite[p. 20]{79.BS}:
\vspace{3mm}

%%%%%%%%%%% L E M M A  B %%%%%%%%%%%%%%%%%%%%%%%%%%%%%%%
\begin{leA}
Let $A_1, \ldots, A_n \in \mathcal{C}_p$ for some $p>0$.
If $A_1, \ldots, A_n$ are positive, then
\begin{equation} \label{ineq:A.p<1}
n^{p-1}\sum_{i=1}^n \|A_i\|_p^p
\le \Norm{\sum_{i=1}^n A_i}_p^p
\le \sum_{i=1}^n \|A_i\|_p^p
\end{equation}
for $0 < p \le 1$ and the reverse inequalities hold for $1 \leq p < \infty$.
\end{leA}
\vspace{3mm}
%%%%%%%%%%%%%%%%%%%%%%%%%%%%%%%%%%%%%%%%%%%%%%%%%%%%%%%%

Let us define a constant $D_{\mathbf{A}}$ for a set of operators $\mathbf{A}=\{A_1,A_2,\ldots,A_n\}$ as follows: \
%%% Def of D_{\mathbf{A}} %%%
\begin{eqnarray*}
D_{\mathbf{A}} := \sum_{i=1}^n \d(A_i) \qquad
\mbox{where} \quad \d(A_i) = \left\{ \begin{array}{@{\,}ll} 1 & (A_i \not= 0) \\ 0 & (A_i = 0) \end{array} \right..
\end{eqnarray*}
%%% Refinement of Lemma A %%%
If there exists $1 \le i \le n$ with $A_i=0$, then Lemma A is refined as follows:
\begin{equation} \label{ineq:DA.p<1}
D_{\mathbf{A}}^{p-1} \sum_{i=1}^n \|A_i\|_p^p
\ \le \ \Norm{\sum_{i=1}^n A_i}_p^p
\ \le \ \sum_{i=1}^n \|A_i\|_p^p
\end{equation}
for $0 < p \le 1$ and the reverse inequalities hold for $1 \leq p < \infty$.\\

We also put $\mathbf{A-B}:=\{A_i-B_j: 1 \le i,j \le n \}$
for sets of operators $\mathbf{A}=\{A_1,A_2,\ldots,A_n\}$ and $\mathbf{B}=\{B_1,B_2,\ldots,B_n\}$.
Then we remark that $0 \le D_{\mathbf{A-B}} \le n^2$.

Now we give our main results that involve three sets of operators.
\vspace{3mm}

%%%%%%%%%%% T H E O R E M  2. 1 %%%%%%%%%%%%%%%%%%%%%%%%
\begin{thm} \label{thm:main}
Let $\mathbf{A}=\{A_1,A_2,\ldots,A_n\},
     \mathbf{B}=\{B_1,B_2,\ldots,B_n\},
     \mathbf{C}=\{C_1,C_2,\ldots,C_n\}
     \subset \mathcal{C}_p$ for some $p>0$.
Then
\begin{equation} \label{ineq:main3.p<2}
\begin{split}
\lefteqn{\sum_{i,j=1}^n \|A_i-A_j\|_p^p +
         \sum_{i,j=1}^n \|B_i-B_j\|_p^p +
         \sum_{i,j=1}^n \|C_i-C_j\|_p^p} \hspace{5mm} \\
&\ge
\left(
  D_{\mathbf{A-B}}^{\frac{p-2}2} \sum_{i,j=1}^n \|A_i-B_j\|_p^p
+ D_{\mathbf{B-C}}^{\frac{p-2}2} \sum_{i,j=1}^n \|B_i-C_j\|_p^p
+ D_{\mathbf{C-A}}^{\frac{p-2}2} \sum_{i,j=1}^n \|C_i-A_j\|_p^p
 \right) \\
&\quad -
\left( \Norm{\sum_{i=1}^n (A_i-B_i)}_p^p
     + \Norm{\sum_{i=1}^n (B_i-C_i)}_p^p
     + \Norm{\sum_{i=1}^n (C_i-A_i)}_p^p \right)
\end{split}
\end{equation}
for $0<p\le 2$ and the reverse inequality holds for $ 2 \le p < \infty$.
\end{thm}

\begin{proof}
We only prove the case when $0 < p \le 2$.
The other case can be proved by a similar argument.\\
We have
\begin{eqnarray*}
&&\hspace{-7mm} \lefteqn{
  \sum_{i,j=1}^n \|A_i-A_j\|_p^p
+ \sum_{i,j=1}^n \|B_i-B_j\|_p^p
+ \sum_{i,j=1}^n \|C_i-C_j\|_p^p} \hspace{5mm} \\
&& + \left( \Norm{\sum_{i=1}^n (A_i-B_i)}_p^p
          + \Norm{\sum_{i=1}^n (B_i-C_i)}_p^p
          + \Norm{\sum_{i=1}^n (C_i-A_i)}_p^p \right) \\
&=&
2\left(\sum_{1\le i<j \le n} \|A_i-A_j\|_p^p
     + \sum_{1\le i<j \le n} \|B_i-B_j\|_p^p
     + \sum_{1\le i<j \le n} \|C_i-C_j\|_p^p\right) \\
&& + \left( \Norm{\sum_{i=1}^n (A_i-B_i)}_p^p
          + \Norm{\sum_{i=1}^n (B_i-C_i)}_p^p
          + \Norm{\sum_{i=1}^n (C_i-A_i)}_p^p \right) \\
&=& 2\left(
  \sum_{1\le i<j \le n} \Norm{|A_i-A_j|^2}_{p/2}^{p/2}
+ \sum_{1\le i<j \le n} \Norm{|B_i-B_j|^2}_{p/2}^{p/2}
+ \sum_{1\le i<j \le n} \Norm{|C_i-C_j|^2}_{p/2}^{p/2}
\right) \\
&& + \left(
\Norm{\left|\sum_{i=1}^n (A_i-B_i)\right|^2}_{p/2}^{p/2}
+
\Norm{\left|\sum_{i=1}^n (B_i-C_i)\right|^2}_{p/2}^{p/2}
+
\Norm{\left|\sum_{i=1}^n (C_i-A_i)\right|^2}_{p/2}^{p/2}
\right) \\
&&\hspace{77mm}
\mbox{(by relation \eqref{mos})} \\
&\ge&
\Norm{\sum_{1\le i<j \le n} |A_i-A_j|^2
      +\sum_{1\le i<j \le n} |B_i-B_j|^2
 +\left| \sum_{i=1}^n(A_i-B_i) \right|^2}_{p/2}^{p/2} \\
&& + \Norm{\sum_{1\le i<j \le n} |B_i-B_j|^2
          +\sum_{1\le i<j \le n} |C_i-C_j|^2
  +\left|\sum_{i=1}^n(B_i-C_i)\right|^2}_{p/2}^{p/2} \\
&& + \Norm{\sum_{1\le i<j \le n} |C_i-C_j|^2
          +\sum_{1\le i<j \le n} |A_i-A_j|^2
  +\left|\sum_{i=1}^n(C_i-A_i)\right|^2}_{p/2}^{p/2} \\
&&\hspace{77mm}
\mbox{(by the second inequality of \eqref{ineq:A.p<1})} \\
&=& \Norm{ \sum_{i,j=1}^n |A_i-B_j|^2 }_{p/2}^{p/2}
  + \Norm{ \sum_{i,j=1}^n |B_i-C_j|^2 }_{p/2}^{p/2}
  + \Norm{ \sum_{i,j=1}^n |C_i-A_j|^2 }_{p/2}^{p/2} \\
&&\hspace{77mm}
\mbox{(by \eqref{oit})} \\
&\ge& D_{\mathbf{A-B}}^{\frac{p}2-1}
      \sum_{i,j=1}^n \Norm{|A_i-B_j|^2}_{p/2}^{p/2}
    + D_{\mathbf{B-C}}^{\frac{p}2-1}
      \sum_{i,j=1}^n \Norm{|B_i-C_j|^2}_{p/2}^{p/2}
    + D_{\mathbf{C-A}}^{\frac{p}2-1}
      \sum_{i,j=1}^n \Norm{|C_i-A_j|^2}_{p/2}^{p/2} \\
&&\hspace{77mm}
\mbox{(by the first inequality of \eqref{ineq:DA.p<1})} \\
&=& D_{\mathbf{A-B}}^{\frac{p-2}2}
    \sum_{i,j=1}^n \|A_i-B_j\|_p^p
  + D_{\mathbf{B-C}}^{\frac{p-2}2}
    \sum_{i,j=1}^n \|B_i-C_j\|_p^p
  + D_{\mathbf{C-A}}^{\frac{p-2}2}
    \sum_{i,j=1}^n \|C_i-A_j\|_p^p.
\end{eqnarray*}
So we have the desired inequality \eqref{ineq:main3.p<2}.
\end{proof}
\vspace{3mm}
%%%%%%%%%%%%%%%%%%%%%%%%%%%%%%%%%%%%%%%%%%%%%%%%%%%%%%%%

The next result can be regarded as a generalization of \eqref{ineq:2-norm}.
\vspace{3mm}

%%%%%%%%%%% P R O P O S I T I O  2. 2 %%%%%%%%%%%%%%%%%%
\begin{prop} \label{pro:main1}
Let $A_1, \ldots, A_n, B_1, \ldots, B_n \in \mathcal{C}_p$ for some $p>0$.
Then
\begin{eqnarray*}
\begin{split}
\lefteqn{\sum_{i,j=1}^n \|A_i-A_j\|_p^p +
\sum_{i,j=1}^n \|B_i-B_j\|_p^p} \hspace{5mm} \\
&\ge 2n^{p-2} \sum_{i,j=1}^n \|A_i-B_j\|_p^p
  - 2\Norm{\sum_{i=1}^n (A_i-B_i)}_p^p
\end{split}
\end{eqnarray*}
for $0 < p \le 2$ and the reverse inequality holds for $ 2 \le p < \infty$.
\end{prop}

\begin{proof}
We only prove the case where $0 < p \le 2$.
Putting $C_i=0$ for $i=1,2,\ldots,n$, we have
$$D_{\mathbf{B-C}}^{\frac{p-2}2} \sum_{i,j=1}^n \|B_i-C_j\|_p^p
= nD_{\mathbf{B}}^{\frac{p-2}2} \sum_{i=1}^n \|B_i\|_p^p, \quad
D_{\mathbf{C-A}}^{\frac{p-2}2} \sum_{i,j=1}^n \|C_i-A_j\|_p^p
= nD_{\mathbf{A}}^{\frac{p-2}2} \sum_{i=1}^n \|A_i\|_p^p.$$
It follows from Theorem \ref{thm:main} that
\begin{eqnarray*}
0 &\le&
  \sum_{i,j=1}^n \|A_i-A_j\|_p^p
+ \sum_{i,j=1}^n \|B_i-B_j\|_p^p
- D_{\mathbf{A-B}}^{\frac{p-2}2} \sum_{i,j=1}^n \|A_i-B_j\|_p^p \\
&& -n\left(
  D_{\mathbf{B}}^{\frac{p-2}2} \sum_{i=1}^n \|B_i\|_p^p
+ D_{\mathbf{A}}^{\frac{p-2}2} \sum_{i=1}^n \|A_j\|_p^p
\right) \\
&& + \left( \Norm{\sum_{i=1}^n (A_i-B_i)}_p^p
          + \Norm{\sum_{i=1}^n B_i}_p^p
          + \Norm{\sum_{i=1}^n A_i}_p^p \right) \\
&=&
  \sum_{i,j=1}^n \|A_i-A_j\|_p^p
+ \sum_{i,j=1}^n \|B_i-B_j\|_p^p
- 2D_{\mathbf{A-B}}^{\frac{p-2}2} \sum_{i,j=1}^n \|A_i-B_j\|_p^p
+ 2\Norm{\sum_{i=1}^n (A_i-B_i)}_p^p \\
&&+ \left( D_{\mathbf{A-B}}^{\frac{p-2}2} \sum_{i,j=1}^n \|A_i-B_j\|_p^p - \Norm{\sum_{i=1}^n (A_i-B_i)}_p^p \right) \\
&&+ \left( \Norm{\sum_{i=1}^n A_i}_p^p
  - nD_{\mathbf{A}}^{\frac{p-2}2} \sum_{i=1}^n \|A_i\|_p^p \right)
+ \left( \Norm{\sum_{i=1}^n B_i}_p^p
  - nD_{\mathbf{B}}^{\frac{p-2}2} \sum_{i=1}^n \|B_i\|_p^p \right).
\end{eqnarray*}
Here the inequality \eqref{ineq:DA.p<1} implies
\begin{eqnarray*}
\lefteqn{D_{\mathbf{A-B}}^{\frac{p-2}2} \sum_{i,j=1}^n \|A_i-B_j\|_p^p - \Norm{\sum_{i=1}^n (A_i-B_i)}_p^p}
\hspace{5mm} \\
&&= D_{\mathbf{A-B}}^{\frac{p}2-1} \sum_{i,j=1}^n \Norm{|A_i-B_j|^2}_{p/2}^{p/2} - \Norm{\sum_{i=1}^n |A_i-B_i|^2}_{p/2}^{p/2} \le 0
\end{eqnarray*}
by relation \eqref{mos}.
Due to $nD_{\mathbf{A}}^{\frac{p-2}2}$ is no less than $1$, we deduce from Lemma A that
$$\Norm{\sum_{i=1}^n A_i}_p^p - nD_{\mathbf{A}}^{\frac{p-2}2} \sum_{i=1}^n \|A_i\|_p^p
\le \Norm{\sum_{i=1}^n A_i}_p^p - \sum_{i=1}^n \|A_i\|_p^p \le 0.$$
Similarly, we have
$\Norm{\sum_{i=1}^n B_i}_p^p \le nD_{\mathbf{B}}^{\frac{p-2}2} \sum_{i=1}^n \|B_i\|_p^p$.
It therefore implies that
\begin{eqnarray*}
\sum_{i,j=1}^n \|A_i-A_j\|_p^p
+ \sum_{i,j=1}^n \|B_i-B_j\|_p^p
&\ge& 2D_{\mathbf{A-B}}^{\frac{p-2}2} \sum_{i,j=1}^n \|A_i-B_j\|_p^p
- 2\Norm{\sum_{i=1}^n (A_i-B_i)}_p^p \\
&\ge& 2n^{p-2} \sum_{i,j=1}^n \|A_i-B_j\|_p^p
- 2\Norm{\sum_{i=1}^n (A_i-B_i)}_p^p \ge 0 \\
&& \hspace{33mm} (\mbox{by $n^2 \ge D_{\mathbf{A-B}} (\ge 0)$}).
\end{eqnarray*}

Thus we obtain the desired inequality.
\end{proof}
\vspace{5mm}
%%%%%%%%%%%%%%%%%%%%%%%%%%%%%%%%%%%%%%%%%%%%%%%%%%%%%%%%

%%%%%%%%%%% C O R O L L A R Y  2. 3 %%%%%%%%%%%%%%%%%%%%
\begin{cor} \label{cor:AB.p-norm}
Let $A_1, \ldots, A_n, B_1, \ldots, B_n \in \mathcal{C}_p$ for some $p>0$ and $\sum_{i=1}^n A_i=\sum_{i=1}^n B_i$.
Then
\begin{equation*}
\sum_{i,j=1}^n \|A_i-A_j\|_p^p +
\sum_{i,j=1}^n \|B_i-B_j\|_p^p
\ge 2n^{p-2} \sum_{i,j=1}^n \|A_i-B_j\|_p^p
\end{equation*}
for $0 < p \le 2$ and the reverse inequality holds for $ 2 \le p < \infty$.
\end{cor}
\vspace{3mm}
%%%%%%%%%%%%%%%%%%%%%%%%%%%%%%%%%%%%%%%%%%%%%%%%%%%%%%%%

Utilizing Corollary \ref{cor:AB.p-norm} with $B_i=0$ ($1 \le i \le n$), we obtain the following result which is a refinement of \cite[Corollary 2.3]{OH.FK.MM}.
\vspace{3mm}

%%%%%%%%%%% C O R O L L A R Y  2. 4 %%%%%%%%%%%%%%%%%%%%
\begin{cor} \label{cor:A.p-norm}
Let $A_1,\ldots,A_n \in \mathcal{C}_p$ for some $p>0$ such that $\sum_{i=1}^n A_i=0$.
Then
\begin{equation*}
\sum_{i,j=1}^n \|A_i-A_j\|_p^p
\ge 2n^{p-1} \sum_{i=1}^n \|A_i\|_p^p
\end{equation*}
for $0 < p \le 2$ and the reverse inequality holds for $ 2 \le p < \infty$.
\end{cor}
\vspace{3mm}
%%%%%%%%%%%%%%%%%%%%%%%%%%%%%%%%%%%%%%%%%%%%%%%%%%%%%%%%

Next, we have the following reverse inequalities of Proposition \ref{pro:main1}:
\vspace{3mm}

%%%%%%%%%%% T H E O R E M  2. 5 %%%%%%%%%%%%%%%%%%%%%%%%
\begin{prop} \label{thm:main2}
Let $A_1, \ldots, A_n, B_1, \ldots, B_n \in \mathcal{C}_p$ for some $p>0$.
Then
\begin{eqnarray*}
\begin{split}
\lefteqn{2\left(n^2-n+1\right)^{\frac{2-p}2}
\sum_{i,j=1}^n \|A_i-B_j\|_p^p
- 2\Norm{\sum_{i=1}^n (A_i-B_i)}_p^p} \hspace{5mm} \\
&\ge \sum_{i,j=1}^n \|A_i-A_j\|_p^p +
    \sum_{i,j=1}^n \|B_i-B_j\|_p^p
\end{split}
\end{eqnarray*}
for $0 < p \le 2$ and the reverse inequality holds for $ 2 \le p < \infty$.
\end{prop}

\begin{proof}
We only consider the case when $0 < p \le 2$.
We have
\begin{eqnarray*}
&&\sum_{1\le i<j\le n} \|A_i-A_j\|_p^p
+ \sum_{1\le i<j\le n} \|B_i-B_j\|_p^p
+ \Norm{\sum_{i=1}^n (A_i-B_i)}_p^p \\
&=& \sum_{1\le i<j \le n} \Norm{|A_i-A_j|^2}_{p/2}^{p/2}
  + \sum_{1\le i<j \le n} \Norm{|B_i-B_j|^2}_{p/2}^{p/2}
  + \Norm{\left|\sum_{i=1}^n (A_i-B_i)\right|^2}_{p/2}^{p/2} \\
  &&\hspace{70mm}
\mbox{(by relation \eqref{mos})} \\
&\le& \left(2\cdot\frac{n^2-n}2+1\right)^{1-\frac{p}2}
\Norm{\sum_{1\le i<j \le n} |A_i-A_j|^2
    + \sum_{1\le i<j \le n} |B_i-B_j|^2
+ \left|\sum_{i=1}^n (A_i-B_i)\right|^2}_{p/2}^{p/2} \\
&&\hspace{70mm}
\mbox{(by the first inequality of \eqref{ineq:A.p<1})} \\
&=& \left(n^2-n+1\right)^{\frac{2-p}2}
    \Norm{\sum_{i,j=1}^n |A_i-B_j|^2}_{p/2}^{p/2}
\hspace{5mm}
\mbox{(by \eqref{oit})} \\
&\le& \left(n^2-n+1\right)^{\frac{2-p}2}
      \sum_{i,j=1}^n \Norm{|A_i-B_j|^2}_{p/2}^{p/2}
\hspace{5mm}
\mbox{(by the second inequality of \eqref{ineq:A.p<1})} \\
&=& \left(n^2-n+1\right)^{\frac{2-p}2} \sum_{i,j=1}^n \|A_i-B_j\|_p^p.
\end{eqnarray*}
So we have the desired inequality.
\end{proof}
\vspace{3mm}
%%%%%%%%%%%%%%%%%%%%%%%%%%%%%%%%%%%%%%%%%%%%%%%%%%%%%%%%

%%%%%%%%%%% R E M A R K  2. 6 %%%%%%%%%%%%%%%%%%%%%%%%%%
\begin{rem}
$\1$
By an easily calculation, we have the inequality
$n^{p-2} < \left(2n^2-n+1\right)^{\frac{2-p}2} < \left(n^{p-2}\right)^{-1}$ for $0 < p \le 2$.

$\2$
The values of $2n^{p-2}$, $2\left(n^2-n+1\right)^{\frac{2-p}2}$ of Propositions \ref{pro:main1} and \ref{thm:main2} are $2$, if $p=2$.
So these results ensured the equality \eqref{ineq:2-norm}.
\end{rem}

Finally we would like to give a problem for further research.
\begin{prob}
What the form of the identity is for the general case of $k$ sets of operators?
\end{prob}

\bigskip
%%%%%%%%%%%%%%%%%%%%%%%%%%%%%%%%%%%%%%%%%%%%%%%%%%%%%%%%
%========================================================================================%
\textbf{Acknowledgment.} The authors sincerely thank the referee for some useful comments. The first author was supported by a grant from Ferdowsi University of Mashhad (No. MP89161MOS).
The third author was supported by Grants-in-Aid for
Scientific Research (No. 20540158), Japan Society for the
Promotion of Science.
%========================================================================================%

%%%%%%%%%%%%%%%%%%%%%%%%%%%%%%%%%%%%%%%%%%%%%%%%%%%%%%%%%           R E F E R E N C E S                           %%%%%%%%%%%%%%%%%%%%%%%%%%%%%%%%%%%%%%%%%%%%%%%%%%%%%%%%

%%%%%%%%%%%%%%%%%%%%%%%%%%%%%%%%%%%%%%%%%%%%%%%%%%%%%%%%

\end{document}